\documentclass[10pt]{article}
\usepackage[latin1]{inputenc}
\usepackage{amsmath,amsfonts,amssymb}
\usepackage{amscd}
\usepackage{amscd,graphicx}

\newcommand{\qed}{{$\diamond$ }}

\begin{document}

\title{Sign conjugacy classes in symmetric groups} 

\author{
J\o rn B. Olsson\\
\small Department of Mathematical Sciences, University of Copenhagen \\[-0.8ex]
\small Universitetsparken 5,DK-2100 Copenhagen \O, Denmark
}

\maketitle

{\small 
{\bf Abstract:} A special type of conjugacy classes in symmetric
groups is studied and used to answer a question about 
odd-degree irreducible characters. 
}

\medskip


\medskip

This work was initiated by a question about associating suitable signs
to odd-degree irreducible characters in the symmetric groups 
$S_n,$ posed by I.M. Isaacs and
G. Navarro. The question is related to their work \cite{IN}.

The positive answer to the question is given below. It is in a sense the
best possible and it involves a special conjugacy class in $S_n.$ 
The has lead the author to a general definition of sign classes in 
finite groups. This general definition is
dicussed briefly in section 1. In section 2 we consider special types of sign
classes in $S_n$ and apply this to the Isaacs-Navarro question in
section 3. The final section contains a general result on sign classes
in $S_n$ and some thoughts about a possible classification of them.

\medskip

{\bf 1.  Sign classes in finite groups}

\medskip

A {\it sign class} in a finite group $G$ is a conjugacy class on which all 
irreducible characters of $G$ take one of the values 0, 1 or
-1. Elements in sign classes are called {\it sign elements}. 

Sign elements of prime order $p$ may occur when you 
have a self-centralizing $p$-Sylow subgroup of order $p$ in $G.$ This
occurs for example for $p=7$ in the simple group $M_{11}$ which also has sign
elements of order 6.  In SL$(2,2^n)$ there is an involution
on which all irreducible characters except the Steinberg character 
take the values 1 or -1. Thus this is a sign element. Non-central
involutions in dihedral groups are also examples of sign elements. 



Column orthogonality for the irreducible characters of $G$
shows that a sign element $s$ gives rise 
to two disjoint multiplicity-free characters $\Theta_s^+$ and $\Theta_s^-$ which
coincide on all conjugacy classes
except the class of $s.$ They are defined as follows 
$$\Theta_s^+=\sum_{\{\chi \in {\rm Irr}(G)| \chi(s)=1\}} \chi ~~{\rm and}~~
\Theta_s^-=\sum_{\{\chi \in {\rm Irr}(G)| \chi(s)=-1\}} \chi $$

An example for symmetric groups is given below. 

Block orthogonality shows that if $p$ is a prime number dividing 
the order of the sign
element $s$ and if you split $\Theta_s^+$ and $\Theta_s^-$ into
components according to the $p$-blocks of characters of $G,$ 
then the values of these components for a given $p$-block still
coincide on all $p$-regular 
elements in $G.$ This has consequences for the decomposition numbers
of $G$ 
at the prime $p.$

\bigskip

{\bf 2. Sign partitions}

\medskip 
 
In this note we are concerned with sign classes in the symmetric groups $S_n.$ 
The irreducible characters of $S_n$ are all integer
valued. Let $\mathcal{P}(n)$ be the set of partitions of $n.$
We write the entries of the character table $X(n)$ of $S_n$  as
$[\lambda](\mu),$ for $\lambda, \mu \in \mathcal{P}(n). $ This is the
value of the irreducible character of $S_n$ labelled by $\lambda$ on
the conjugacy class labelled by $\mu.$  

We call  $\mu \in \mathcal{P}(n)$ a {\it sign partition} if the  
the corresponding conjugacy class is a sign class, i.e. if 
$[\lambda](\mu) \in \{0,1,-1\}$ for all $\lambda \in \mathcal{P}(n). $
The {\it support} of a sign partition
$\mu$ is defined as 
$${\rm supp}(\mu)=\{\lambda \in \mathcal{P}(n)~|~[\lambda](\mu) \neq 0\}$$

For example  the
Murnaghan-Nakayama formula (\cite{JK}, 2.4.7 or \cite{J} 21.1 shows 
that $(n)$ is always a sign partition. Indeed 
$[\lambda](n) \neq 0$ if and only if $\lambda=(n-k, 1^k)$ is a hook
partition and then  $[\lambda](n)=(-1)^k.$
Using column orthogonality for irreducible characters this has as a 
consequence that the generalized character
$$\Theta_{(n)}=\sum_{k=0}^{n-1} (-1)^k [n\!-\!k, 1^k]$$ 
takes the value 0 everywhere except on the class $(n)$ where is has
value $n.$ 

For an arbitrary sign partition $\mu$
$$\Theta_{\mu}= \sum_{\lambda \in {\rm supp}(\mu)} [\lambda](\mu) [\lambda]$$ 
is a  generalized character vanishing outside the conjugacy class
of $\mu$ and it is the difference between disjoint multiplicity-free characters 
$\Theta_{\mu}^+$ and $\Theta_{\mu}^-$. (See section 1.)

\medskip

Below is a list of all sign partitions for $n=2,...,10:$

\medskip

\noindent $n=2:~   (2), (1^2)$
 
\noindent $n=3: ~  (3), (2,1)$

\noindent $n=4:~   (4), (3, 1), (2, 1^2)$

\noindent $n=5: ~  (5), (4, 1),  (3, 2), (3, 1^2)$

\noindent  $n=6: ~  (6), (5, 1), (4, 2), (4, 1^2), (3, 2, 1)$

\noindent  $n=7: ~  (7), (6, 1), (5, 2), (5, 1^2),  (4,3), (4,2,1),
(3, 2, 1^2) $

\noindent $n=8:  ~ (8), (7, 1), (6, 2), (6, 1^2),  (5, 3), (5, 2, 1), 
(4 ~3~ 1)$

\noindent $n=9: ~  (9), (8, 1), (7, 2), (7, 1^2), (6,3), (6, 2, 1),  
 (5, 4), (5, 3 , 1), (5, 2 ,1^2)$

\noindent $n=10: (10), (9, 1), (8, 2), (8, 1^2), (7, 3), (7, 2, 1), 
(6, 4), (6, 3, 1), $

$~~~~~~~~(6, 2 ,1^2), (5,4,1), (4,3,2,1)$

\medskip

The sign partition (4, 2) of 6 yields two characters
of degree 20
$$\Theta_{(4,2)}^+=[6]+[4,2]+[2^2,1^2]+[1^6] $$ 
$$\Theta_{(4~2)}^-=[5,1]+[3^2]+[2^3]+[2,1^4] $$ 
coinciding everywhere except on the class (4,2) where they differ by a
sign.

\medskip

An important class of sign partitions are the {\it unique
  path}-partitions (for short {\it up}-partitions). They are described as
follows. 
 If  $\mu=(a_1,a_2,...,a_k)$ and $\lambda$ are  partitions of $n$,
 then a $\mu$-path in $\lambda$ is a sequence $\lambda=\lambda_0,
 \lambda_1,..., \lambda_k=(0),$ of partitions, where for $i=1...k$ $\lambda_i$ is
 obtained by removing an $a_i$-hook in $\lambda_{i-1}.$ 
Then we call $\mu$ is an {\it up}-partition for $\lambda$ if the number of 
$\mu$-paths in $\lambda$ is at most 1.
We call $\mu$ is an {\it up}-partition if it is a {\it up}-partition
for all partitions $\lambda$ of $n.$

\medskip

\noindent {\bf Proposition 1:} {\it A up-partition is also a sign partition.}

\medskip

\noindent {\it Proof:} This follows immidiately by repeated use of the
Murnaghan-Nakayama formula. 
If there is no $\mu$-path for $\lambda,$ then
$[\lambda](\mu)=0.$ Otherwise  $[\lambda](\mu)=(-1)^k$, where $k$ is
the sum of the leg lengths of the hooks involved in the unique
$\mu$-path for $\lambda.$ \qed

\medskip

\noindent {\bf Remarks:} 1. If $\mu=(a_1,a_2,...,a_k)$ is an {\it up}-partition
with $a_k=2,$ then also 
$\mu'=(a_1,a_2,...,a_{k-1}, 1^2)$ is an {\it up}-partition. 

2. If $\mu=(a_1,a_2,...,a_k)$ is an {\it up}-partition
with $k \ge 2,$ then also 
$\mu^*=(a_2,...,a_k)$ is an {\it up}-partition. 
Indeed, if a partition $\lambda^*$ of $n-a_1$ has two or more
$\mu^*$-paths then a partition of $n$ obtained be adding an $a_1$-hook
to $\lambda^*$ has two or more $\mu$-paths.

3. The partition $(3,2,1)$ is s sign partition, but not a 
{\it up}-partition, since there are two  $(3,2,1)$-paths in 
the partition $(3,2,1).$ 
Also $(4,3,2,1)$ is a sign partition, but not a 
{\it up}-partition, since there are two  $(4,3,2,1)$-paths in 
the partition $(7,2,1).$  

\medskip

\noindent {\bf Proposition 2:} {\it Let $m>n.$ If 
$\mu^*=(a_1,a_2,...,a_k)$ is a partition of $n,$  and  $\mu=(m, a_1,a_2,...,a_k)$
then $\mu^*$ is a sign partition (respectively a {\it up}-partition) of
$n$ if and 
only if $\mu$ is a sign partition (respectively a {\it
  up}-partition) of $m+n.$  }

\medskip

\noindent {\it Proof:} Let $\lambda$ be a partition of $m+n.$ 
Since $2m>m+n$ $\lambda$ cannot contain more
than at most one hook of length $m,$ e.g.  by 2.7.40 in \cite{JK}.
This clearly implies that $\mu^*$ is a {\it up}-partition if and only
of  $\mu$ is a {\it up}-partition. 
If $\lambda$ has no hook of length $m,$ then $[\lambda](\mu)=0.$
If $\lambda$ has a hook of length $m,$ then remove the unique hook of
that length to get the partition $\lambda_1.$ Then
$[\lambda](\mu)=\pm  [\lambda_1](\mu^*).$ 
If $\mu^*$ is a sign
partition we get that  $[\lambda_1](\mu^*) \in \{0,1,-1\}$ and thus  
$[\lambda](\mu) \in \{0,1,-1\}.$ This shows that if $\mu^*$ is a sign
partition then $\mu$ is a sign partition. If $\mu$ is a sign partition
and if $\lambda_1\in \mathcal{P}(n),$ then add a hook of length $m$ to
$\lambda_1$ to get a partition $\lambda.$ Since by assumption
$[\lambda](\mu) \in \{0,1,-1\},$ the same is true for
$[\lambda_1](\mu^*).$   \qed

\medskip

It is an intersting question whether it is possible to recognize from the parts of
$\mu,$  whether or not $\mu$ is an {\it up}-partition or a sign partition. The final
section of this paper contains results related to this question. 

However the above proposition suggests the following definition of a
class of sign partitions, given in terms of its parts. 

If $\mu=(a_1,a_2,...,a_k)$ is a partition we call it {\it strongly
  decreasing } (for short a ${\it sd}$-partition) if  
we have $a_i > a_{i+1}+...+a_k$ for $i=1,...,k-1.$ 

\medskip

\noindent {\bf Remarks:} 1. Obviously, if  $\mu=(a_1,a_2,...,a_k)$
is an {\it sd}-partition with $k \ge 2$ then $\mu^*=(a_2,...,a_k)$
is also an {\it sd}-partition.

2. The partition $(3,1^2)$ is an {\it up}-partition, but not an {\it sd}-partition.

\medskip

\noindent {\bf Proposition 3:} {\it An sd-partition is a
  up-partition and thus also sign partition.}

\medskip

\noindent {\it Proof:} That an {\it sd}-partition is a
  {\it up}-partition is proved by repeated use of Proposition 2. \qed

\medskip

\noindent {\bf Remark:} The {\it sd}-partitions are closely related to
the so-called ``non-squashing'' partitions. A partition
$\mu=(a_1,a_2,...,a_k)$ is called non-squashing if  $a_i \ge
a_{i+1}+...+a_k$  for all $i=1,...,k-1.$ It is known that 
that the number non-quashing partitions of $n$ equals the 
binary  partitions of $n,$ i.e. the number of  partitions of
$n$ into parts which are powers of 2. (\cite{HS}, \cite{SS}).
Let $s(n)$ denote the number of {\it sd}-partitions
of $n.$ Put $s(0)=1.$ Ordering the set of {\it sd}-partitions according to their largest
part shows that 
$$s(n)=\sum_{i=0}^{\lfloor (n-1)/2 \rfloor}s(i).$$
Thus for all $k \ge 1$  we have $s(2k-1)=s(2k).$ Putting 
$t(k)=2s(2k)=s(2k-1)+s(2k)$ it can be shown that $t(k)$ is then equal
to the number
of binary partitions of $2k.$ 

\medskip

\noindent {\bf Proposition 4:} {\it If  $\mu=(a_1,a_2,...,a_k)$ is an
sign partition of $n$ then the number of irreducible characters
$\lambda$ with $[\lambda](\mu) \neq 0$ is  $z_{\mu},$  the order of the
centralizer of an element of type $\mu$ in $S_n.$ In particular, for an
sd-partition  $z_{\mu}=a_1a_2...a_k.$}

\medskip

\noindent {\it Proof:} Since the non-zero values of irreducible
characters on $\mu$ are 1 or -1 this follows from column orthogonality. 
\qed

\bigskip

{\bf 3. The Isaacs-Navarro question.}
 
\medskip

Some background for this may be found in \cite{IN}.

\smallskip 

{\bf Question:} (Isaacs-Navarro) 
 Let $P$ be 2-Sylow subgroup of $S_n$ and  $Irr_{2'}(S_n)$ be the set
of odd degree irreducible characters of $S_n$. Does there exist
signs $e_{\chi}$ for $\chi \in Irr_{2'}(S_n)$ such that the character 
  $$\Theta= \sum_{\chi \in Irr_{2'}(S_n)}e_{\chi}\chi$$
satisfies that 
  $$(i)~~~~\Theta(x) {\rm ~is ~divisible ~by} ~ |P/P'| ~ {\rm for ~ all} ~x \in S_n.$$
and 
 $$ (ii)~~~~\Theta(x)=0 ~ {\rm for ~all} ~ x \in S_n ~{\rm of~ odd~ order}?   $$

This is answered positively by

\medskip

\noindent {\bf Theorem 5:} {\it Write
  $n=2^{r_1}+2^{r_2}+...+2^{r_t}$, where $r_1>r_2>...>r_t \ge 0.$
Then $\mu=(2^{r_1},2^{r_2},...,2^{r_t}) $ is a sd-partition 
with support {\rm supp}$(\mu)=Irr_{2'}(S_n).$ Moreover $\Theta_{\mu}$ satisfies the 
conditions (i) and (ii) above. Indeed   $\Theta_{\mu}$ vanishes
everywhere except on $\mu$ where it takes the value  $|P/P'|.$ }

 \medskip

\noindent {\it Proof:} Clearly $\mu$ is an {\it sd}-partition and thus
a sign partition, which implies that  $\Theta_{\mu}$ vanishes
everywhere except on $\mu$ where it takes the value
$z_{\mu}=2^{r_1+r_2+...+r_t}.$ This is the cardinality of supp($\mu)$
(Proposition 4).
If $C_i$ is the iterated wreath product
of $i$ copies of the cyclic group of order 2 then $C_i/C_i'$ is an
elementary abelian group of order $2^i. $ Since $P \simeq C_{r_1}
\times C_{r_2} \times ... \times  C_{r_t}$ we get
$|P/P'|=2^{r_1+r_2+...+r_t}.$ We need then only the fact that
supp$(\mu)= Irr_{2'}(S_n).$ 
By \cite{MNO}, Theorem 4.1, supp$(\mu) \subseteq Irr_{2'}(S_n).$ 
(Since we here know that non-zero values on $\mu$ are $\pm 1,$ this also 
follows from a general character theoretic result, \cite{F},(6.4))
On the other hand $|Irr_{2'}(S_n)|=2^{r_1+r_2+...+r_t}$ by \cite{M}, Corollary (1.3), so
that the supp($\mu$) cannot be properly contained in  
$Irr_{2'}(S_n).$ \qed

\medskip

\noindent {\bf Remark:} The results from  \cite{M}, \cite{MNO} quoted in
the above proof are formulated for arbitrary primes. However Theorem 5
does not have an analogue for odd primes.

\medskip

\noindent {\bf Example:} In $SL(2,2^n)$ the 2-Sylow subgroup is self
centralizing. It has a unique conjugacy class of involutions and $2^n+1$ irreducible  
characters, all of which (with the exception of the Steinberg character)
have odd degrees. The involutions are sign element, so that
$\Theta_t,$ $t$ involution, vanishes on all elements of odd order. the
value on $t$ is $2^n.$ Thus this is another example of the existence
of signs for odd degree irreducible characters such that the signed
sum satisfy the conditions mentioned above.

\bigskip

{\bf 4. Repeated parts in a sign partition}

\medskip

We want to show that repeated parts are very rare in sign
partitions. Indeed only the part 1 may be repeated. 

\medskip

\noindent {\bf Lemma 6:} {\it A sign partition $\mu$ cannot have its smallest
  part repeated except for the part 1, which may be repeated once.}
 
\medskip

\noindent {\it Proof:} Suppose that 1 is repeated $m \ge 2$ times in
$\mu$ then
$[n-1,1](\mu)=[m-1,1](1^m)=m-1. $ Thus $m=2.$
If $b>1$ is the smallest part, repeated   $m \ge 2$ times then 
$[n-b,b](\mu)=m.$ \qed

\medskip

\noindent {\bf Theorem 7:} {\it A sign partition cannot have
  repeated parts except for the part 1, which may be repeated once.}
 
\medskip

\noindent {\it Proof:} We are going to assume that $a$ is the smallest
repeated part $>1$ in the partition $\mu$ and that the multiplicity of $a$ in
$\mu$ is $m\ge 2.$ We want to determine a partition $\lambda$
satisfying that {\it all hook lengths outside the first row  are $\le a$ } and 
in addition  $|[\lambda](\mu)| \ge m.$

Divide the parts of $\mu$ into 

$a_1 \ge...\ge a_{i-1}$ (all greater than $a$) (sum $s$, say)

$a_i,...,a_{i+m-1}$ ($m$ parts all equal to $a$)  

$a_{i+m}>...>a_k.$ (all parts smaller than $a$) (sum $t,$ say)
(However we allow $a_k-1=a_k=1.$ )

We let $\mu^*=(a_{i+m},...,a_k).$

By Lemma 6 we may assume that $t>0.$  An easy analysis shows that we 
may assume {\bf $a \ge 4.$} 
(To do this we just have to show that partitions on the form $$
(2^m,1),(2^m,1^2), (3^m,2,1), (3^m,2,1^2),(3^m,1), (3^m,1^2), m\ge 2$$
are not sign partitions. For example $[2m-1,1^2](2^m,1)=[2m,1^2](2^m,1^2)=-m.$)

First we notice that we need only consider the case that $s=0.$ 
Indeed, if $\lambda'$ is a partition of $n-s$ satisfying that all hook
lengths outside the first row  are $\le a$ and that  $|[\lambda'](\mu')| \ge m,$ 
where $\mu'=(a_i,...,a_k)$ and $\lambda$ is obtained by adding $s$ to
the largest part of $\lambda'$ then MN shows that
$[\lambda](\mu)=[\lambda'](\mu')$ and we are done. 
(Here and in the the following MN refers to the Murnaghan-Nakayama formula.)
Thus we may assume that $a=a_1$ is the only repeated part of $\mu,$
apart possibly from 1.

We have then $n=ma+t.$ Let for $0\le i \le m$ $\mu_i$ be $\mu$ with
$i$ parts $a$ removed. Thus $\mu_0=\mu$ and $\mu_m=\mu^*.$

Now $(n-a,1^a)$ has only two hooks of length $a$ so MN shows  
$[n-a,1^a](\mu)=(-1)^{a-1}[n-a](\mu_1)+[n-2a,1^a](\mu_1)=(-1)^{a-1}+[n-2a,1^a](\mu_1).$

Inductively we get $[n-a,1^a](\mu)=(m-1)(-1)^{a-1}+[t,1^a](\mu_{m-1}).$

If {\bf $t \le a$} then $[t,1^a]$ has only one hook of length $a$ and we get
 $[t,1^a](\mu_{m-1})=(-1)^{a-1}[t](\mu_m)=(-1)^{a-1}$ and thus  
 $[n-a,1^a](\mu)=m(-1)^{a-1}$. Thus $[n-a,1^a]$ may be chosen as the
 desired $\lambda.$ 

We may assume {\bf $a < t.$} 

Consider the case $t < 2a$ so that $t-a < a.$ There are exactly $a$
partitions of $t$ obtained by adding an $a$-hook to the
partition $(t-a).$ Suppose that $\kappa_i$ is obtained by adding a
hook with leg length $i$ to  $(t-a).$ 

Since  $t < 2a$ each $\kappa_i$ has only one hook of length $a$ 
(eg.  by 2.7.40 in \cite{JK}).
Removing it we get $(t-a).$ Note that $\kappa_0=(t).$
By Theorem 21.7 in \cite{J} the generalized character $\sum_{i=0}^{a-1}(-1)^i
\kappa_i$ takes the value 0 on $\mu^*,$ since $\mu^*$ has no part
divisible by $a.$ Choose an $j>0$ such that 
$(-1)^j[\kappa_j](\mu^*) \geq 0.$ (Clearly, the
$(-1)^j[\kappa_j](\mu^*)$ cannot all be $<0$, since the contribution 
from $[a]$ is equal to 1 and $a \geq 4.$) Put 
$\lambda^*=\kappa_j$
so that 
$$(-1)^j [\lambda^*](\mu^*) \geq 0.$$
Let $\lambda$ be obtained from $\lambda^*$ by adding $ma$ to its
largest part. Thus the largest part of $\lambda$ is at least $n-a$ so
that trivially all hook lengths outside the first row are $\le a.$ 
We claim that $|[\lambda](\mu)| \ge m.$ 

Let for $0 \le i \le m$ $\lambda_i$ be obtained by subtracting $ia$ from the
largest part of $\lambda,$ so that $\lambda_0=\lambda$ and
$\lambda_m=\lambda^*.$ Let $\mu_i$ be as above.

By MN we have 
$$[\lambda_i](\mu_{i})=[\lambda_{i+1}](\mu_{i+1})+(-1)^{j}$$
for $0\le 1 <m.$ Thus
$$[\lambda](\mu)=[\lambda_{1}](\mu_{1})+(-1)^{j}$$
$$=[\lambda_{2}](\mu_{2})+2(-1)^{j}$$
and so on. This shows   
$$[\lambda](\mu)=[\lambda^*](\mu^*)+m(-1)^{j}.$$
Thus 
$$[\lambda](\mu)=[\lambda^*](\mu^*)+m(-1)^{j}=(-1)^j((-1)^j[\lambda^*](\mu^*)+m). $$
This has absolute value $\geq m,$ so that $\mu$ is not a sign class.

\smallskip

A similar argument may be used in the case  $t \ge 2a.$ Then $t-a \ge
a$ and it is possible to add a $a$-hook to the partition $(t-a)$ in
$a+1$ ways. Putting an $a$-hook with leg length $i$ below $t-a$ gives
you $a$ partitions  $\kappa_i, i=0,...,i-1.$ In addition we have the
partition $(t)$. Using again Theorem 21.7 in \cite{J} we see that 
the generalized character $\sum_{i=0}^{a-1}(-1)^i
\kappa_i$ takes the value -1 on $\mu^*.$ It is possible to
choose an $j\geq 0$ such that $(-1)^j [\kappa_j](\mu^*) \geq 0.$ Otherwise we
would have  $-1=\sum_{i=0}^{a-1}(-1)^i[\kappa_i](\mu^*) \leq -a. $ We then
proceed as in the previous case. \qed

\medskip

\noindent {\bf Corollary 8:} {\it If $\mu$ is a sign partition, then the
  centralizer of elements of cycle type $\mu$ is abelian. In short:
  Centralizers of sign elements in $S_n$ are abelian.}

\medskip

\noindent {\bf Remark.} G. Navarro has kindly pointed out that there
exists a group of order 32 containing a sign element with a 
non-abelian centralizer.

\medskip

\noindent {\bf Corollary 9:} {\it Suppose that 
 $n=2^{r_1}+2^{r_2}+...+2^{r_t}$, where $r_1>r_2>...>r_t \ge 0.$ 
The sign classes of 2-elements in $S_n$ have for $n$ odd (i.e. $r_t=0$) 
cycle type $(2^{r_1},2^{r_2},...,2^{r_t}).$ If $n=4k+2$ (i.e. $r_t=1)$ we have in addition 
$(2^{r_1},2^{r_2},...,2^{r_{t-1}},1^2).$ If $n=8k+4$  (i.e. $r_t=2)$ we have in
  addition $(2^{r_1},2^{r_2},...,2^{r_{t-1}},2,1^2).$}

\medskip

\noindent {\it Proof:} If a sign class of an 2-element in $S_n$ does
not have the type $(2^{r_1},2^{r_2},...,2^{r_t}),$ then by 
Theorem 7, the part 1 has to be repeated twice. We have seen that
$(1^2)$ and $(2,1^2)$ are sign partitions. Therefore Proposition 2
shows that the other two cycle types listed in the corollary are indeed 
cycle types for sign classes. For these values of $n$ there can be no
more sign classes. If $n=8k,$ (i.e. $r_t \ge 3$) then the possibility
that $2^{r_t}$ is replaced by $2^{r_{t-1}},2^{r_{t-2}},...,2,1^2$ is
  excluded by  Proposition 2 and the fact that $(4,2,1^2)$ is not a
  sign partition. \qed

\medskip

Finally we formulate a conjecture about which partitions are sign
partitions. It seems that sign partitions are close to being {\it
  sd}-partitions. 

We fix the following notation:  
 $\mu^*=(a_1,a_2,...,a_r)$ for some $r \ge 2$ is a partition
of $t$ and $\mu=(a,a_1,...,a_r)$
where $a>a_1.$ 

Then $\mu$ is called {\it exceptional} if $a \le t$ and both $\mu$ and $\mu^*$
are sign partitions and in addition the partitions
$(a_i,a_{i+1},...,a_r)$ are all sign partitions. 

If we can determine the exceptional partitions, then we also know all the
sign partitions. However there exist infinite series of exceptional
partitions. Indeed it can be shown that the following partitions are
exceptional:

$\bullet$ $(a,a-1,1)$ for $a \ge 2.$

$\bullet$ $(a,a-1,2,1)$ for $a \ge 4.$

$\bullet$ $(a,a-1,3,1)$ for $a \ge 5.$

\medskip

The author suspects strongly that these are the only infinite series
of exceptional partitions and would
like to state the following conjecture.

\medskip

\noindent {\bf Conjecture:}  Let $\mu=(a_1,a_2,...,a_k)$ be a partition. Then
$\mu$ is a sign partition if and only if one of the following
conditions hold: 

{\rm (1)} $\mu$ is an {\it sd}-partition, i.e.  $a_i>a_{i+1}+...+a_k$ for $i=1,...,k-1.$

{\rm (2)}  $a_i>a_{i+1}+...+a_k$ for $i=1,...,k-2$ and in addition $a_{k-1}=a_k=1.$

{\rm (3)}  $a_i>a_{i+1}+...+a_k$  for $i=1,...,k-3$ and  in addition
$(a_{k-2},a_{k-1},a_k)=(a,a-1,1)$ for some $a \ge 2.$
 
{\rm (4)} $a_i>a_{i+1}+...+a_k$ for  $i=1,...,k-4$ and  in addition
$(a_{k-3},a_{k-2},a_{k-1},a_k)$ is one of the following

$\bullet$ $(a,a-1,2,1)$ for some $a \ge 4$

$\bullet$ $(a,a-1,3,1)$ for some $a \ge 5$

$\bullet$ (3,2,1,1)

$\bullet$ (5,3,2,1).

\medskip

We hope to return to this conjecture in a later paper. Its
verification would also easily imply a classification 
of {\it up}-partitions.

\medskip

{\bf Acknowledgments:} The author thanks G. Navarro for  
the question, which initiated this work and C. Bessenrodt for 
some discussions. Part of this work was done during the authors visit 
to the Mathematical Sciences Research Institute (MSRI) in April-May 2008.

\vspace{3ex}

\end{document}